\documentclass[10pt]{amsart}  
\usepackage{amsmath,amsthm,amssymb,cite}
\begin{document}
\title{STABILITY THEOREMS FOR GROUP ACTIONS ON UNIFORM SPACES} 
\author{Pramod Das, Tarun Das} 
\begin{abstract}
We extend the notions of topological stability, shadowing and persistence from homeomorphisms to finitely generated group
actions on uniform spaces and prove that an expansive action with either shadowing or persistence is topologically stable. Using the concept of null set of a Borel measure $\mu$, we introduce the notions of $\mu$-expansivity, $\mu$-topological stability, $\mu$-shadowing and $\mu$-persistence for finitely generated group actions on uniform spaces and show that a $\mu$-expansive action with either $\mu$-shadowing or $\mu$-persistence is $\mu$-topologically stable.     
\end{abstract}        
\maketitle

\textbf{Mathematics Subject Classifications (2010):} 37C75, 37C85, 37C50, 49J53, 54H20    
\medskip

\textbf{Keywords and Phrases:} Expansive Homeomorphisms, Expansive Measures, Shadowing, Persistence, Topological Stability
\medskip

$\textbf{1. INTRODUCTION}$
\medskip

Topological stability is a fundamental notion of a dynamical system which guarantees that the qualitative behaviour of trajectories remains unaffected by continuous small perturbations. In Lyapunov stability one considers perturbations of initial conditions in a fixed system. Unlike Lyapunov stability, in topological stability one considers perturbations of the system itself. Walter's stability theorem applied to Anosov diffeomorphisms of compact smooth manifolds is one of the finest results in differentiable dynamics.    
\medskip

\textbf{Theorem 1.1}[\cite{W},\textit{Theorem 1}] Let $M$ be a compact smooth Manifold and let $f:M\rightarrow M$ be a diffeomorphism. If $f$ is Anosov then it is topologically stable. 
\medskip

The following result can be regarded as an extension of \textit{Theorem 1.1}, to homeomorphisms on compact metric spaces because of the fact that Anosov diffeomorphisms are expansive and has shadowing \cite{R}.  
\medskip

\textbf{Theorem 1.2} [\cite{AH}, Theorem 2.4.5] Let $X$ be a compact metric space and let $f:X\rightarrow X$ be a homeomorphism. If $f$ is expansive and has shadowing then it is topologically stable.
\medskip

In \cite{CCL}, authors have improved \textit{Theorem 1.2} by weakening the hypothesis of the theorem. In particular, one can observe that persistence is weaker than shadowing. 
\medskip

\textbf{Theorem 1.3} [\cite{CCL},\textit{Theorem 2}] Let $X$ be a compact metric space and let $f:X\rightarrow X$ be a homeomorphism. If $f$ is expansive and persistent then it is topologically stable. 
\medskip

Second author of the present paper with others generalized \textit{Theorem 1.2}, to homeomorphisms on non-compact, non-metrizable spaces. 
\medskip

\textbf{Theorem 1.4} [\cite{DLRW},\textit{Theorem 22}] Let $X$ be a first countable, locally compact, para-compact, Hausdorff uniform space and let $f:X\rightarrow X$ be a homeomorphism. If $f$ is expansive and has shadowing then it is topologically stable. 
\medskip

In \cite{CL} and \cite{ALL}, authors have extended \textit{Theorem 1.2} and \textit{Theorem 1.3} respectively, to finitely generated group actions on compact metric spaces.
\par\noindent\rule{\textwidth}{0.4pt} 
\begin{tiny}
Department of Mathematics, Faculty of Mathematical Sciences, University of Delhi, New Delhi-110007. 
\\
Email Address: tarukd@gmail.com (Tarun Das), pramod.math.ju@gmail.com (Pramod Das)
\end{tiny} 

\textbf{Theorem 1.5} [\cite{CL},\textit{Theorem 2.8}] Let $X$ be a compact metric space and let $\Phi$ be an action of a finitely generated group on $X$. If $\Phi$ is expansive and has shadowing then it is topologically stable. 
\medskip

\textbf{Theorem 1.6} [\cite{ALL},\textit{Theorem A}] Let $X$ be a compact metric space and let $\Phi$ be an action of a finitely generated group on $X$. If $\Phi$ is expansive and persistent then it is topologically stable. 
\medskip

Our first aim is to extend \textit{Theorem 1.5} and \textit{Theorem 1.6} to finitely generated group actions on non-compact, non-metrizable spaces.    
\medskip

\textbf{Theorem 1.7} Let $X$ be a first countable, locally compact, para-compact, Hausdorff uniform space and let $\Phi$ be an action of a finitely generated group on $X$. If $\Phi$ is expansive and has either shadowing or persistence then it is topologically stable.       
\medskip

In recent past, the notion of expansivity has been extended to include Borel measures \cite{M}. On the other hand, Lee and Morales have introduced \cite{LM} the notions of shadowing and topological stability for Borel measures on compact metric spaces. They have proved that expansive measures with shadowing are topologically stable. 
\medskip

\textbf{Theorem 1.8} [\cite{LM},\textit{Theorem 3.1}] Let $X$ be a compact metric space and let $f:X\rightarrow X$ be a homeomorphism. If a measure $\mu$ is expansive and has shadowing with respect to $\Phi$ then it is topologically stable.
\medskip

Very recently, authors \cite{DKY} have extended this theorem to finitely generated group actions on compact metric spaces. 
\medskip

\textbf{Theorem 1.9} [\cite{DKY}, \textit{Theorem 3.2}] Let $X$ be a compact metric space and let $\Phi$ be an action of a finitely generated group on $X$. If $\Phi$ is $\mu$-expansive and has $\mu$-shadowing then it is $\mu$-topologically stable.  
\medskip

Our second aim is to extend these results to finitely generated group actions on non-compact, non-metrizable spaces. 
\medskip

\textbf{Theorem 1.10} Let $\Phi$ be an action of a finitely generated group on a first countable, locally compact, para-compact, Hausdorff uniform space $X$ and let $\mu$ be a Borel measure on $X$. If $\Phi$ is $\mu$-expansive and has either $\mu$-shadowing or $\mu$-persistence then it is $\mu$-topologically stable.        
\medskip

By following the name of the very first result among these (the differentiable case), it is reasonable to name all the above mentioned results as \textit{Stability Theorems}.           
\medskip

We organize this paper as follows. In section 2, we define the notions of topological stability, $\mu$-topological stability, shadowing, $\mu$-shadowing, expansivity, $\mu$-expansivity, persistence, $\mu$-persistence for finitely generated group actions on uniform spaces. In section 3, we prove \textit{Theorem 1.7} and \textit{Theorem 1.9} and discuss examples. In section 4, we discuss general properties of the introduced notions.      
\medskip

\textbf{2. PRELIMINARIES AND DEFINITIONS} 
\medskip

Let $X$ be a non-empty set. Then, the diagonal of $X\times X$ is given by $\Delta(X)=\lbrace (x,x)|x\in X\rbrace$. For a subset $R$ of $X\times X$, we define $R^{-1}=\lbrace (y,x)|(x,y)\in R\rbrace$. We say that $R$ is symmetric if $R=R^{-1}$. For two subsets $U$ and $V$ of $X\times X$, we define their composition as $U\circ V=\lbrace(x,y)\in X\times X |$ there is $z\in X$ satisfying $(x,z)\in U$ and $(z,y)\in V\rbrace$. In this paper, we assume that the phase space of a dynamical system is a uniform space $(X,\mathcal{U})$, where $\mathcal{U}$ is a collection of subsets of $X\times X$ satisfying the following properties \cite{K}:  
\begin{enumerate}
\item[(i)] Every $D\in\mathcal{U}$ contains $\Delta(X)$.
\item[(ii)] If $D\in\mathcal{U}$ and $E\supset D$, then $E\in\mathcal{U}$.
\item[(iii)] If $D,D'\in\mathcal{U}$, then $D\cap D'\in\mathcal{U}$.
\item[(iv)] If $D\in\mathcal{U}$, then $D^{-1}\in\mathcal{U}$. 
\item[(v)] For every $D\in\mathcal{U}$ there is a symmetric $D'\in\mathcal{U}$ such that $D'\circ D'\subset D$.
\end{enumerate}

The members of $\mathcal{U}$ are called entourages of $X$. If $(X,\mathcal{U})$ is a uniform space, then we can generate a topology on $X$ by characterizing that a subset $Y\subset X$ is open if and only if there is an entourage $U$ of $X$ such that for each $x\in Y$ the cross section $U[x]=\lbrace y\in X\mid (x,y)\in U\rbrace$ is contained in $Y$. An entourage $M$ is said to be proper if for any compact set $K\subset X$, $M[K]=\bigcup_{x\in K}M[x]$ is compact. An entourage $U$ is said to be wide if there is a compact subset $K\subset X$ such that $U\cup (K\times X)=X\times X$. In other words, $U$ is wide if there is a compact set $K\subset X$ such that for any $x\notin K$, we have $U[x]=X$. A point $x\in X$ is called an atom if $\mu(\lbrace x\rbrace)>0$. A measure $\mu$ on $X$ is said to be non-atomic if it has no atom. 
\medskip

In this paper, $G$ denotes a finitely generated group and $X$ a uniform space with uniformity $\mathcal{U}$. Any generating set $S$ of the group $G$ will be considered finite symmetric (for any $s\in S$, $s^{-1}\in S$).     
\medskip

A map $\Phi:G\times X\rightarrow X$ is said to be an action of $G$ on $X$ if the following holds:
\begin{enumerate}
\item[(i)] The map $\Phi_g=\Phi(g,.)$ is a uniform equivalence for any $g\in G$.
\item[(ii)] $\Phi_e(x)=x$ for all $x\in X$, where $e$ is the identity element of the group $G$. 
\item[(iii)] $\Phi_{g_1g_2}(x)=\Phi_{g_1}(\Phi_{g_2}(x))$ for all $x\in X$ and $g_1,g_2\in G$. 
\end{enumerate}

We denote $Act(G,X)$, the set of actions of a finitely generated group $G$ on the uniform space $X$.  
\medskip

$\textbf{Definition 2.1}$ Let $X$ be a uniform space with uniformity $\mathcal{U}$. Then, 
\\
(i) An action $\Phi\in Act(G,X)$ is called expansive if there exists closed entourage $D\in\mathcal{U}$ such that for any distinct $x,y\in X$, there exists $g\in G$ such that $(\Phi_g(x),\Phi_g(y))\notin D$. Such $D$ is called an expansive entourage for $\Phi$.
\\
(ii) An action $\Phi\in Act(G,X)$ is called $\mu$-expansive for some non-atomic measure $\mu$, if there exists closed entourage $D\in\mathcal{U}$ such that $\mu(\Gamma_D(x))=0$ for all $x\in X$, where $\Gamma_D(x)=\lbrace y\in X\mid (\Phi_g(x),\Phi_g(y))\in D$ for all $g\in G\rbrace$. Such $D$ is called a $\mu$-expansive entourage for $\Phi$. 
\medskip

\textbf{Example 2.2} (i) Let $(\mathbb{R},d)$ be the real line with usual topology generated by the usual metric $d$. Let us give $\mathbb{R}^k$ the topology generated by the metric 
\\
$D((x_1,x_2,...,x_k),(y_1,y_2,...,y_k))=$max$\lbrace d(x_1,y_1),d(x_2,y_2),..,$ $d(x_k,y_k)\rbrace$. 
\\
Let $\Phi:\mathbb{Z}^k\times \mathbb{R}^k\rightarrow \mathbb{R}^k$ be given by $\Phi_{(n_1,n_2,...,n_k)}(x_1,x_2,...,x_k)=(2^{n_1}x_1,2^{n_2}x_2,...,2^{n_k}x_k)$. Then one can check that $\Phi$ is expansive. 
\\
(ii) Let $X$ be a uniform space and let $f:X\rightarrow X$ be a $\mu$-expansive uniform equivalence. Then, the actions $\Phi:\mathbb{Z}^k\times X\rightarrow X$ by $\Phi_{(n_1,n_2,...,n_k)}(x)=f^{n_1+n_2+...+n_k}x$ and $\Psi:\mathbb{Z}^k\times X\rightarrow X$ given by $\Psi_{(n_1,n_2)}(x)=f^{n_1}(x)$ are $\mu$-expansive.      
\medskip

\textbf{Definition 2.3} (i) Let $\Phi\in Act(G,X)$ and let $D\in\mathcal{U}$ be given. Then, $\lbrace x_g\rbrace_{g\in G}$ is said to be a $D$-pseudo orbit (with respect to the generating set $S$) for $\Phi$ if $(x_{sg},\Phi_s(x_g))\in D$ for all $s\in S$, $g\in G$. On the other hand, $\lbrace x_g\rbrace_{g\in G}$ is said to be $D$-shadowed by some point $x\in X$ if $(x_g,\Phi_g(x))\in D$ for all $g\in G$. A $D$-pseudo orbit $\lbrace x_g\rbrace_{g\in G}$ is said to be through a measurable set $B$ if $x_e\in B$, where $e$ is the identity element of the group $G$.  
\\
(ii) Let $\mathcal{P}(X)$ be the power set of $X$ and let $H:X\rightarrow \mathcal{P}(X)$ is a set valued map of $X$. We define the domain of $H$ by $Dom(H)=\lbrace x\in X\mid H(x)\neq\phi\rbrace$. $H$ is said to be compact valued if $H(x)$ is compact for each $x\in X$. We write $(Id,H)\in D$ for some $D\in\mathcal{U}$ if $H(x)\subset D[x]$ for each $x\in X$. $H$ is called upper semi-continuous if for every $x\in Dom(H)$ and every open neighborhood $O$ of $H(x)$ there exists $D\in\mathcal{U}$ such that $H(y)\subset O$ for all $x,y\in X$ with $(x,y)\in D$.     
\medskip

\textbf{Definition 2.4} Let $X$ be a uniform space with uniformity $\mathcal{U}$. Then, 
\\
(i) An action $\Phi\in Act(G,X)$ is said to have shadowing (with respect to the generating set $S$ of $G$) if for every $E\in\mathcal{U}$ there exists $D\in\mathcal{U}$ such that every $D$-pseudo orbit is $E$-shadowed by some $x\in X$. 
\\
(ii) An action $\Phi\in Act(G,X)$ is said to be topologically stable (with respect to the generating set $S$ of $G$) if for every $E\in\mathcal{U}$ there exists $D\in\mathcal{U}$ such that if $\Psi\in Act(G,X)$ is another action satisfying $(\Psi_s(x),\Phi_s(x))\in D$ for all $x\in X$, $s\in S$ then there exists a continuous map $f:X\rightarrow X$ such that $\Phi_g\circ f=f\circ\Psi_g$ for all $g\in G$ and $(x,f(x))\in E$ for all $x\in X$.
\\
(iii) An action $\Phi\in Act(G,X)$ is said to have $\mu$-shadowing for some $\mu\in M(X)$ (with respect to the generating set $S$ of $G$), if for every $E\in\mathcal{U}$ there exists $D\in\mathcal{U}$ and a measurable set $B\subset X$ with $\mu(X\setminus B)=0$ such that every $D$-pseudo orbit through $B$ is $E$-shadowed by some point $x\in X$.  
\\
(iv) An action $\Phi\in Act(G,X)$ is said to be $\mu$-topologically stable for some $\mu\in M(X)$ (with respect to the generating set $S$ of $G$), if for every $E\in\mathcal{U}$ there exists $D\in\mathcal{U}$ such that if $\Psi\in Act(G,X)$ is another action satisfying $(\Psi_s(x),\Phi_s(x))\in D$ for all $x\in X$, $s\in S$ then there is an upper semi-continuous compact valued map $H$ of $X$ with measurable domain such that the following conditions hold: (a) $\mu(X\setminus Dom(H))=0$, (b) $\mu\circ H=0$, (c) $(Id,H)\in E$, (d) $\Phi_g\circ H=H\circ\Psi_g$ for all $g\in G$. 
\\
(v) An action $\Phi\in Act(G,X)$ is said to be persistent (with respect to the generating set $S$ of $G$) if for every $E\in\mathcal{U}$ there exists $D\in\mathcal{U}$ such that if $\Psi\in Act(G,X)$ is another action satisfying $(\Psi_s(x),\Phi_s(x))\in D$ for all $x\in X$, $s\in S$ then for every $x\in X$ there exists $y\in X$ such that $(\Psi_g(x),\Phi_g(y))\in E$ for all $g\in G$.  
\\
(vi) An action $\Phi\in Act(G,X)$ is said to be $\mu$-persistent for some $\mu\in M(X)$ (with respect to the generating set $S$ of $G$) if for every $E\in\mathcal{U}$ there exists $D\in\mathcal{U}$ and a measurable set $B\subset X$ with $\mu(X\setminus B)=0$ such that if $\Psi\in Act(G,X)$ is another action satisfying $(\Psi_s(x),\Phi_s(x))\in D$ for all $x\in X$, $s\in S$ then for each $x\in B$ there exists $y\in X$ such that $(\Psi_g(x),\Phi_g(y))\in E$ for all $g\in G$.    
\medskip

\textbf{Lemma 2.5} Let $S$ and $T$ be two distinct generating sets of $G$. Then, for any $\Phi\in Act(G,X)$  
\begin{enumerate}
\item[(i)] $\Phi$ has shadowing with respect to $S$ if and only if it has shadowing with respect to $T$.
\item[(ii)] $\Phi$ is topologically stable with respect $S$ if and only if it is topologically stable with respect $T$.
\item[(iii)] $\Phi$ is persistent with respect to $S$ if and only if it is persistent with respect to $T$.
\item[(iv)] $\Phi$ has $\mu$-shadowing with respect to $S$ if and only if it has $\mu$-shadowing with respect to $T$.
\item[(v)] $\Phi$ is $\mu$-topologically stable with respect to $S$ if and only if is $\mu$-topologically stable with respect to $T$.
\item[(vi)] $\Phi$ is $\mu$-persistent with respect to $S$ if and only if it is $\mu$-persistent with respect to $T$.    
\end{enumerate}

$\textbf{Proof.}$ (i) Suppose $\Phi$ has shadowing with respect to $S$. Let $E\in\mathcal{U}$ be given and $D'\in\mathcal{U}$ be given for $E$ by the shadowing of $\Phi$. It suffices to show that there exists $D\in\mathcal{U}$ such that every $D$-pseudo orbit with respect to $T$ is $D'$-pseudo orbit with respect to $S$. Put $m:=$max$_{s\in S} l_T(s)$, where $l_T$ is the word length metric on $G$ induced by $T$. Choose $D_1\in\mathcal{U}$ such that $D_1^m\subset D'$. Since $\Phi_g$ is uniform equivalence for each $g\in G$ and $S$, $T$ finite, there exists $D\subset D_1$ such that $(\Phi_h(x),\Phi_h(y))\in D_1$ for all $x,y\in X$ with $(x,y)\in D$ and $h\in G$ with $l_T(h)\leq m$. Any $s\in S$ can be written as $s=t_1t_2t_3...t_{l(s)}$, where $l(s)=l_T(s)\leq m$, $t_i\in T$, $i=1,2,3,...,l(s)$. Let $\lbrace y_g\rbrace_{g\in G}$ be a $D$-pseudo orbit for $\Phi$ with respect to $T$, i.e; $(y_{tg},\Phi_t(y_g))\in D$ for all $t\in T$, $g\in G$. 
\medskip

Claim: $(y_{sg},\Phi_s(y_g))=(y_{t_1...t_{l(s)}g},\Phi_{t_1...t_{l(s)}}(y_g))\in D'$ for all $s\in S$, $g\in G$. Observe that 
\begin{center}
$(y_{t_1...t_{l(s)}g},\Phi_{t_1}(y_{t_2...t_{l(s)}g}))\in D\subset D_1$, $(\Phi_{t_1}(y_{t_2...t_{l(s)}g}),\Phi_{t_1}(\Phi_{t_2}(y_{t_3...t_{l(s)}g}))\in D_1$,..., 
\newline
$(\Phi_{t_1...t_{l(s)-1}}(y_{l(s)g}),\Phi_{t_1...t_{l(s)-1}}(\Phi_{l(s)}(y_g)))\in D_1$. 
\end{center}

Then, the claim holds by the definition of composition of $m$-number of entourages and the fact that $D_1^m\subset D'$. 
\medskip

(ii) Suppose $\Phi$ is $S$-topologically stable. Let $E\in\mathcal{U}$ be given and $D'\in\mathcal{U}$ be given for $E$ by the topological stability of $\Phi$. It suffices to show that there exists $D\in\mathcal{U}$ such that for any $\Psi\in Act(G,X)$ if $(\Psi_t(x),\Phi_t(x))\in D$ for all $t\in T$, $x\in X$ then $(\Psi_s(x),\Phi_s(x))\in D'$ for all $s\in S$, $x\in X$. Put $m:=$max$_{s\in S} l_T(s)$, where $l_T$ is the word length metric on $G$ induced by $T$. Choose $D_1\in\mathcal{U}$ such that $D_1^m\subset D'$. Since $\Phi_g$ is uniform equivalence for each $g\in G$ and $S$, $T$ are finite, there exists $D\in\mathcal{U}$ such that $(\Phi_h(x),\Phi_h(y))\in D_1$ for all $x,y\in X$ with $(x,y)\in D$ and $h\in G$ with $l_T(h)\leq m$. Any $s\in S$ can be written as $s=t_1t_2t_3...t_{l(s)}$, where $l(s)=l_T(s)\leq m$, $t_i\in T$, $i=1,2,3,...,l(s)$.
\\

Claim: $(\Psi_s(x),\Phi_s(x))=(\Psi_{t_1...t_{l(s)}}(x),\Phi_{t_1...t_{l(s)}}(x))\in D_1^m\subset D'$ for all $s\in S$, $x\in X$.  
\begin{center}
Observe that $(\Psi_{t_1...t_{l(s)-1}}\Psi_{t_{l(s)}}(x),\Psi_{t_1...t_{l(s)-1}}\Phi_{t_{l(s)}}(x))\in D_1$, $(\Psi_{t_1...t_{l(s)-2}}\Psi_{t_{l(s)-1}}\Phi_{t_{l(s)}}(x)),\Psi_{t_1...t_{l(s)-2}}\Phi_{t_{l(s)-1}}\Phi_{t_{l(s)}}(x)))\in D_1...,$
\\
$(\Psi_{t_1}\Phi_{t_2...t_{l(s)}}(x),\Phi_{t_1...t_{l(s)}}(x))\in D_1$.
\end{center}
Then, the claim holds by the definition of composition of $m$ number of entourages and the fact that $D_1^m\subset D'$. 
\medskip

The proofs for (iii), (iv), (v) and (vi) follow similarly as those of (i) and (ii).    
\medskip

In view of the above lemma, the following are natural definitions of shadowing, topological stability, $\mu$-shadowing, $\mu$-topological stability, persistence and $\mu$-persistence for $\Phi\in Act(G,X)$.            
\medskip

$\textbf{Definition 2.6}$ (i) An action $\Phi\in Act (G,X)$ is said to have shadowing if it has shadowing with respect to some generating set $S$ of $G$.
\\
(ii) An action $\Phi\in Act(G,X)$ is said to be topologically stable if it is topologically stable with respect to some generating set $S$ of $G$.
\\
(iii) An action $\Phi\in Act(G,X)$ is said to be persistent if it is persistent with respect to some generating set $S$ of $G$. 
\\
(iv) An action $\Phi\in Act(G,X)$ is said to have $\mu$-shadowing if it has $\mu$-shadowing with respect to some generating set $S$ of $G$.  
\\
(v) An action $\Phi\in Act(G,X)$ is said to be $\mu$-topologically stable if it is $\mu$-topologically stable with respect to some generating set $S$ of $G$. 
\\
(vi) An action $\Phi\in Act(G,X)$ is said to be $\mu$-persistent if it is $\mu$-persistent with respect to some generating set $S$ of $G$.
\medskip

\textbf{Remark 2.7} Following \cite{DD}, we call the above mentioned notions as topological version of the corresponding metric notions on compact metric spaces (see \cite{CL} for expansivity, shadowing and topological stability, \cite{ALL} for persistence, \cite{DKY} for $\mu$-expansivity, $\mu$-shadowing and $\mu$-topological stability. If $(X,d)$ is compact metric space, then the topological notions coincide with the corresponding metric notions. For deeper understanding of this fact one may refer to \textit{Remark 3.2}\cite{DD}.   
\medskip

Observe that shadowing implies $\mu$-shadowing, topological stability implies persistence. We now show that $\mu$-topological stability implies $\mu$-persistence.  
\medskip

\textbf{Proposition 2.9} If an action $\Phi\in Act(G,X)$ is $\mu$-topologically stable, then it is $\mu$-persistent.  
\medskip

$\textbf{Proof.}$ Let $E\in\mathcal{U}$ be symmetric and let $D\in\mathcal{U}$ be given for $E$ by the $\mu$-topological stability of $\Phi$. Let $\Psi\in Act(G,X)$ be another action such that $(\Psi_s(x),\Phi_s(x))\in D$ for all $x\in X$, $s\in S$. Let $H$ be the set-valued map as in the definition of $\mu$-topological stability of $\Phi$. If $B=Dom(H)$, then $\mu(X\setminus B)=0$ and for every $x\in B$, there exists $y\in X$ such that $y\in H(x)$ which implies $\Phi_g(y)\in\Phi_g(H(x))$ for all $g\in G$. This gives $\Phi_g(y)\in H(\Psi_g(x))\subset E[\Psi_g(x)]$ for all $g\in G$ which implies $(\Psi_g(x),\Phi_g(y))\in E$ for all $g\in G$. That means $\Phi$ is $\mu$-persistent.  
\medskip

The following lemma was proved in \cite{DLRW}. For sake of completeness, we give the proof here as well. 
\medskip 

\textbf{Lemma 2.10} Let $X$ be a first countable, locally compact, para-compact, Hausdorff uniform space and let $\mu$ be a non-atomic measure on $X$. If $\Phi\in Act(G,X)$ is expansive ($\mu$-expansive), then there is a proper expansivity ($\mu$-expansivity) entourage for $\Phi$.  
\medskip
 
\textbf{Proof.} Let $E\in\mathcal{U}$ be an expansivity ($\mu$-expansivity) entourage for $\Phi$. Since $X$ is locally compact, each $x\in X$ has an open, relatively compact neighborhood $U_x$. Since $X$ is para-compact, the open cover $\lbrace U_x\rbrace_{x\in X}$ has a closed (and hence compact) locally finite refinement ${V_{\alpha}}$ [\cite{K}, Chapter 5, \textit{Theorem 28}]. Then, $A=E\cap (\bigcup_{\alpha} (V_{\alpha}\times V_{\alpha}))$ is a proper expansivity ($\mu$-expansivity) entourage.  
\pagebreak

\textbf{3. STABILITY THEOREMS}  
\medskip

In this section, our first aim is to prove the following theorem which extends \textit{Theorem 2.8} \cite{CL}.            
\medskip

\textbf{Theorem 3.1} Let $X$ be a first countable, locally compact, para-compact, Hausdorff uniform space. If $\Phi\in Act(G,X)$ is expansive and has shadowing, then it is topologically stable. 
\medskip

In fact, we prove a stronger result in the next theorem and get the above theorem as a direct consequence because of the fact that actions with shadowing are persistent. Further, we know that topologically stable actions are persistent. So, it is reasonable to find sufficient conditions for a persistent action to be topologically stable. So, the following result which extends $\textit{Theorem 3.2}$\cite{CCL}, is important in its own right.     
\medskip

\textbf{Theorem 3.2} Let $X$ be a first countable, locally compact, para-compact, Hausdorff uniform space. If $\Phi\in Act(G,X)$ is expansive and persistent, then it is topologically stable.   
\medskip

To prove the above theorem we need the following lemmas of which $\textit{Lemma 3.4}$ was proved in \cite{DLRW} and the proof of $\textit{Lemma 3.3}$ is left for the reader as an easy exercise. 
\medskip

\textbf{Lemma 3.3} Let $E\in\mathcal{U}$ be such that $E^2$ is expansive entourage for $\Phi\in Act(G,X)$. Let $D\in\mathcal{U}$ be given for $E$ by persistence of $\Phi$. If $\Psi\in Act(G,X)$ is another action with $(\Psi_s(x),\Phi_s(x))\in D$ for all $x\in X$, $s\in S$ then for every $x\in X$ there is unique $y\in X$ such that $(\Psi_g(x),\Phi_g(y))\in E$ for all $g\in G$.     
\medskip

\textbf{Lemma 3.4} Let $X$ be a first countable, locally compact, para-compact, Hausdorff uniform space. A map $f:X\rightarrow X$ is continuous if for any wide $U\in\mathcal{U}$, there is $V\in\mathcal{U}$ such that $(f\times f)(V)\subset U$. 
\medskip

\textbf{Lemma 3.5} Let $X$ be a first countable, locally compact, para-compact, Hausdorff uniform space and let $\Phi$ be an action of $G$ on $X$ with a proper expansive entourage $A$. For any non-empty finite set $F\subset G$, define $V_F(A)=\lbrace (x,y)\in X\times X\mid (\Phi_g(x),\Phi_g(y)\in A$ for all $g\in F\rbrace$. Then, for any wide $U\in\mathcal{U}$ there exists a non-empty finite set $F\subset G$ such that $V_F(A)\subset U$.  
\medskip

$\textbf{Proof.}$ By contradiction, suppose there exists a wide $U\in\mathcal{U}$ such that for each non-empty finite set $F\subset G$ there exists $(x_F,y_F)\in V_F(A)\cap ((X\times X)\setminus U)$. Choose a sequence of non-empty finite subsets $F_1\subset F_2\subset F_3\subset...$ of $G$ such that $G=\bigcup_{n\geq 1} F_n$. Let $L=\lbrace (x_n,y_n)\in V_{F_n}(A)\cap ((X\times X)\setminus U)\mid n\geq 1\rbrace$. Since $U\in\mathcal{U}$ is wide, then there exists a compact set $K\subset X$ such that $U\cup (K\times X)=X\times X$. Thus, $L\subset ((X\times X)\setminus U)\cap A\subset K\times A[K]$, which is a compact subset of $X\times X$. So, there exists a sequence $(x_{n_k},y_{n_k})\in L$ converging to $(x,y)\notin\Delta(X).$ On the other hand, we have $(\Phi_g(x),\Phi_g(y))\in A$ for all $g\in G$, which leads to a contradiction to the fact that $A$ is an expansive entourage for $\Phi$.      
\medskip

\textbf{Proof of Theorem 3.2} Let $B\in\mathcal{U}$ be given and let $A$ be a proper expansive entourage for $\Phi$. Choose symmetric entourage $E$ such that $E^3\subset A\cap B$. Since $\Phi$ is persistent, there exists $D\in\mathcal{U}$ such that if $\Psi\in Act(G,X)$ is another action with $(\Psi_s(x),\Phi_s(x))\in D$ for all $x\in X$, $s\in S$, by $\textit{Lemma 3.3}$ for every $x\in X$, there is unique $y\in X$ such that $(\Psi_g(x),\Phi_g(y))\in E$ for all $g\in G$. Define a map $f:X\rightarrow X$ as given by $f(x)=y$. Observe that $(x,f(x))\in E\subset E^3\subset B$ for all $x\in X$.      
\medskip

We now show that $f(\Psi_h(x))=\Phi_h(f(x))$ for all $x\in X$, $h\in G$. We have $(\Psi_{gh}(x),\Phi_g(f(\Psi_h(x)))=(\Psi_g(\Psi_h(x)),\Phi_g(f(\Psi_h(x)))\in E$ for all $g, h\in G$. On the other hand, $(\Psi_{gh}(x),\Phi_g(\Phi_h(f(x))))=(\Psi_{gh}(x),\Phi_{gh}(f(x)))\in E$ for all $g,h\in G$. Since $E$ is symmetric, $(\Phi_g(f(\Psi_h(x)),\Phi_g(\Phi_h(f(x)))\in E^2\subset A$ for all $g,h\in G$. Since $A$ is expansive entourage, we have $f(\Psi_h(x))=\Phi_h(f(x))$ for all $x\in X$, $h\in G$.            
\medskip

Finally, we prove that $f:X\rightarrow X$ is continuous by using $\textit{Lemma 3.4}$. Let $U$ be a wide entourage. Then, by $\textit{Lemma 4.5}$ there exists non-empty finite set $F\subset G$ such that $V_F(A)\subset U$. Let $W=\bigcap_{g\in F}(\Psi_g\times \Psi_g)(E)$ and we show that $(f\times f)(W)\subset U$. Let $(x,y)\in W$. Then, for any $g\in F$, $(\Psi_g(x),\Phi_g(f(x)))\in E$, $(\Psi_g(x),\Psi_g(y))\in E$ and $(\Psi_g(y),\Phi_g(f(y)))\in E$. Since $E$ is symmetric, $(\Phi_g(f(x)),\Phi_g(f(y))=(f(\Psi_g(x)),f(\Psi_g(y)))\in E^3\subset A$ which implies $(f\times f)(W)\subset V_F(A)\subset U$. This completes the proof. 
\medskip

\textbf{Example 3.6} Let $G=<a,b\mid ba=a^2b>$ be an infinite solvable group and $S=\lbrace a,b,a^{-1},b^{-1}\rbrace$ be a generating set for $G$. Let $\Phi:G\times\mathbb{R}^n\rightarrow\mathbb{R}^n$ be generated by the maps $\Phi_a((x_1,x_2,...,x_n))=(x_1,x_2,...,x_n)$, $\Phi_{a^{-1}}((x_1,x_2,...,x_n))=(x_1,x_2,...,x_n)$, $\Phi_b((x_1,x_2,...,x_n))=(mx_1,mx_2,...,mx_n)$ and $\Phi_{b^{-1}}((x_1,x_2,...,x_n))=(\frac{1}{m}x_1,\frac{1}{m}x_2,...,\frac{1}{m}x_n)$, where $m>1$. One can follow the same steps as in \textit{Proposition 2.2} \cite{B} to show that $\Phi$ has shadowing. It is also easy to check that $\Phi$ is expansive. So, by \textit{Theorem 3.1} $\Phi$ is topologically stable.  
\medskip

\textbf{Example 3.7} Let $\mathbb{R}$ be the real line and let $f:\mathbb{R}\rightarrow \mathbb{R}$ be the homeomorphism given by $f(x)=2x$. Then, the action $\Phi\in Act(\mathbb{Z},X)$ given by $\Phi_n(x)=f^n(x)$ is expansive and has shadowing. Therefore, by $\textit{Theorem 3.1}$, $\Phi$ is topologically stable.     
\medskip

Our next aim is to prove the following measurable stability theorem which extends \textit{Theorem 3.1} \cite{LM}.       
\medskip

\textbf{Theorem 3.8} Let $X$ be first countable, locally compact, para-compact, Hausdorff uniform space and let $\mu$ be a non-atomic measure on $X$. If $\Phi\in Act(G,X)$ is $\mu$-expansive and has $\mu$-shadowing, then it is $\mu$-topologically stable.    
\medskip

We prove a stronger result in the following theorem and obtain the above theorem as a direct consequence. From $\textit{Proposition 2.9}$, we see that every $\mu$-topologically stable action is $\mu$-persistent. Thus, the following theorem provides a sufficient condition for a $\mu$-persistent action to be $\mu$-topologically stable.    
\medskip

\textbf{Theorem 3.9} Let $X$ be a first countable, locally compact, para-compact, Hausdorff uniform space and let $\mu$ be a non-atomic measure on $X$. If $\Phi\in Act(G,X)$ is $\mu$-persistent and $\mu$-expansive, then it is $\mu$-topologically stable. 
\medskip

$\textbf{Proof.}$ Let $E\in\mathcal{U}$ be given and let $A$ be a $\mu$-expansive entourage for $\Phi$. Further, let $F\in\mathcal{U}$ be a closed symmetric entourage such that $F^2\subset A\cap E$. Then, by \textit{Lemma 2.10} there is a proper entourage $E'\in\mathcal{U}$ such that $E'^2\subset F$. Let $D\in\mathcal{U}$ and $B\subset X$ be given for $E'$ by $\mu$-persistence of $\Phi$. Fix $\Psi\in Act(G,X)$ such that $(\Psi_s(x),\Phi_s(x))\in D$ for all $x\in X$, $s\in S$. Define a compact-valued map $H:X\rightarrow \mathcal{P}(X)$ by $H(x)=\bigcap_{g\in G}\Phi_{g^{-1}}(E'[\Psi_g(x)])$. 
\medskip

Let us first prove that $Dom(H)$ is measurable. Let $\lbrace x_k\rbrace_{k\in\mathbb{N}}$ be a sequence in $Dom(H)$ such that $x_k\rightarrow x$ as $k\to\infty$. Then, for each $k\in\mathbb{N}$, we can choose $y_k$ such that $(\Psi_g(x_k),\Phi_g(y_k))\in E'$ for all $g\in G$. In particular, $y_k\in E'[x_k]$ for all $k\in\mathbb{N}$. By local compactness and Hausdorff property of $X$, there exists a compact neighborhood $K$ of $x$. Since $x_k\to x$ as $k\to\infty$, there is $N\in\mathbb{N}$ such that $x_k\in K$ for all $k\geq N$. Therefore, $y_k\in E'[K]$ for all $k\geq N$ and since $E'[K]$ is compact, $y_k$ converges to a point, say $y$ in $X$. Taking $k\to\infty$, we have that $(\Psi_g(x),\Phi_g(y))\in E'$ for all $g\in G$ which implies that $y\in H(x)$. Thus, $H(x)\neq\phi$ and hence, $x\in Dom(H)$. By the first countability of $X$ we have $Dom(H)$ is closed and so is measurable.         
\medskip

Next, we prove that $\mu(X\setminus Dom(H))=0$. Since $(\Psi_s(x),\Phi_s(x))\in D$ for all $x\in X$, $s\in S$, for fixed $x\in B$, there exists $y\in X$ such that $(\Psi_g(x),\Phi_g(y))\in E'$ for all $g\in G$. This shows that $H(x)\neq\phi$ for each $x\in B$ and hence $B\subset Dom(H)$. Therefore, $\mu(X\setminus Dom(H))\leq\mu(X\setminus B)=0$.  
\medskip

We now prove that $H$ is upper semi-continuous. Fix $x\in Dom(H)$ and an open neighborhood $O$ of $H(x)$. Since $G$ is finitely generated, there are non-empty subsets $F_1\subset F_2\subset F_3\subset...$ such that $G=\bigcup_{n\geq 1} F_n$. Define $H(y)=\bigcap_{m\geq 1}H_m(y)$, where $H_m(y)=\bigcap_{g\in F_m}\Phi_{g^{-1}}(E'[\Psi_g(y)])$. Observe that each $H_m(y)$ is compact and $H_{m+1}(y)\subset H_m(y)$ for all $m\geq 1$. Taking $y=x$, we obtain that $H_m(x)\subset O$ for all $m\geq 1$. We assert that there exists $D'\in\mathcal{U}$ such that $H_m(y)\subset O$ for all $y\in X$ with $(x,y)\in D'$. If not, then there exists a sequence $y_k\rightarrow x$ and $z_k\in H_m(y_k)\setminus O$ for all $k\geq 1$. By the local compactness and Hausdorff property of $X$, $z_k$ converges to a point, say $z$ and observe that $z\notin O$. But $z_k\in H_m(y_k)$ implies that $(\Psi_g(y_k),\Phi_g(z_k))\in E'$ for all $k\geq 1$ and $g\in F_m$. Letting $k\to\infty$, we have $(\Psi_g(x),\Phi_g(z))\in E'$ for all $g\in F_m$. So, $z\in H_m(x)\subset O$, a contradiction. By the assertion $H(y)\subset H_m(y)\subset O$, whenever $(x,y)\in D'$. Thus, $H$ is upper semi-continuous.   
\medskip

We now prove that $\mu\circ H=0$, i.e; $\mu(H(x))=0$ for all $x\in X$. Fix $x\in X$ and let $y\in H(x)$. If $z\in H(x)$, then $(\Psi_g(x),\Phi_g(z))\in E'\subset F$ for all $g\in G$. Since $y\in H(x)$, we have $(\Psi_g(x),\Phi_g(y))\in E'\subset F$ for all $g\in G$. Since $F$ is symmetric, we get $(\Phi_g(z),\Phi_g(y))\in F^2\subset A$ for all $g\in G$. Since $A$ is symmetric, $H(x)\subset\Gamma_A(y)$ and since $A$ is a $\mu$-expansive entourage for $\Phi$, $\mu(H(x))\leq \mu(\Gamma_A(y))=0$ for all $x\in X$. Also, observe that $H(x)\subset E'[x]$ and since $E'\subset E$, we have that $(Id,H)\in E$.  
\medskip

Finally, we prove that $\Phi_h\circ H=H\circ\Psi_h$ for all $h\in G$. If $x\in Dom(H)$, then $H(x)\neq\phi$ and hence $\Phi_h(H(x))=\Phi_h(\bigcap_{g\in G}\Phi_{g^{-1}}(E'[\Psi_g(x)]))=\bigcap_{g\in G} \Phi_{hg^{-1}}(E'[\Psi_g(x)])=\bigcap_{g\in G}\Phi_{g^{-1}}(E'[\Psi_{gh}(x)])=\bigcap_{g\in G}\Phi_{g^{-1}}(E'[\Psi_g(\Psi_h(x))])=H(\Psi_h(x))$ for all $h\in G$.  On the other hand, if $x\notin Dom(H)$, then $\Psi_h(x)\notin Dom(H)$ and hence $\Phi_h(H(x))=\phi=H(\Psi_h(x))$. This completes the proof.     
\medskip

\textbf{Example 3.10} Let $f:X\rightarrow X$ be an $N$-expansive homeomorphism with shadowing \cite{CC} on a compact metric space. Then, the action $\Phi\in Act(\mathbb{Z},X)$ given by $\Phi_n(x)=f^n(x)$ is $\mu$-expansive with shadowing. So by $\textit{Theorem 3.8}$, $\Phi$ is $\mu$-topologically stable.                                  
\medskip

\textbf{Theorem 3.11} If $\Phi\in Act(G,X)$ is topologically stable, then it is $\mu$-topologically stable for any non-atomic measure $\mu$ on $X$ (if it exists).     
\medskip

\textbf{Proof.} Suppose $\Phi\in Act(G,X)$ is topologically stable and $\mu$ be a non-atomic measure on $X$. We want to show that $\Phi$ is $\mu$-topologically stable. Let $E$ be given symmetric entourage and let $D\in\mathcal{U}$ be given for $E$ by the topological stability of $\Phi$. Fix $\Psi\in Act(G,X)$ be another action such that $(\Psi_s(x),\Phi_s(x))\in D$ for all $x\in X$, $s\in S$. Then, there exists a continuous map $f:X\rightarrow X$ such that $\Phi_g\circ f=f\circ\Psi_g$ for all $g\in G$. Define the compact valued map $H:X\rightarrow \mathcal{P}(X)$ given by $H(x)=\lbrace f(x)\rbrace$ for all $x\in X$. Then, $H$ is upper semi-continuous because of the fact that $f$ is continuous. Observe that $Dom(H)=X$ and hence, $\mu(X\setminus Dom(H))=0$. Since $\mu$ is non-atomic $\mu(H(x))=\mu(\lbrace f(x)\rbrace)=0$ for all $x\in X$ proving $\mu\circ H=0$. Since $(x,f(x))\in E$ for all $x\in X$, $(x,H(x))\in E$ for all $x\in X$ and hence $(Id,H)\in E$. Further, since $\Phi_g\circ f=f\circ \Psi_g$ for all $g\in G$, we must have $\Phi_g\circ H=H\circ\Psi_g$ for all $g\in G$. This completes the proof.  
\medskip

\textbf{Corollary 3.12} A complete separable metric space admitting topologically stable action which is not $\mu$-topologically stable is at most countable.    
\medskip

\textbf{Proof.} This follows from $\textit{Theorem 3.11}$ and $\textit{Theorem 8.1, P. 53}$ in \cite{P}. 
\medskip

\textbf{4. GENERAL PROPERTIES OF INTRODUCED NOTIONS}  
\medskip

An action $\Phi\in Act(G,X)$ is said to be uniformly conjugate to another action $\Psi\in Act(G,Y)$ if there is a uniform equivalence $h:X\rightarrow Y$ such that $h\circ\Phi_g=\Psi_g\circ h$ for all $g\in G$. 
\medskip

$\textbf{Proposition 4.1}$ Let $\Phi\in Act(G,X)$ and $\Psi\in Act(G,Y)$ be uniformly conjugate. Then,
\begin{enumerate}
\item[(i)] $\Phi$ has shadowing if and only if $\Psi$ has shadowing.   
\item[(ii)] $\Phi$ is expansive if and only if $\Psi$ is expansive. 
\item[(iii)] $\Phi$ is persistent if and only if $\Psi$ is persistent. 
\item[(iv)] $\Phi$ is topologically stable if and only if $\Psi$ is topologically stable.     
\end{enumerate} 
\medskip

\textbf{Proof.} We prove (i) and (iii) only. The proofs of (ii) and (iv) are left for the reader as easy exercises.  
\medskip

(i) Suppose that $\Phi\in Act(G,X)$ has shadowing. Let $h:X\rightarrow Y$ be a uniform equivalence satisfying $h\circ\Phi_g=\Psi_g\circ h$ for all $g\in G$. Let $E\subset Y\times Y$ be given entourage and let the entourage $D\subset X\times X$ be given for $E$ by the uniform continuity of $h$, i.e; $(h(x),h(y))\in E$ whenever $(x,y)\in D$. Let the entourage $A\subset X\times X$ be given for $D$ by the shadowing of $\Phi$ and since $h^{-1}:Y\rightarrow X$ is uniform equivalence, there exists another entourage $B\subset Y\times Y$ such that $(h^{-1}(x),h^{-1}(y))\in A$ whenever $(x,y)\in B$.  
\medskip

Let $\lbrace y_g\rbrace_{g\in G}$ be $B$-pseudo orbit for $\Psi$, i.e; $(y_{sg},\Psi_s(y_g))\in B$ for all $g\in G$, $s\in S$. This implies $(h^{-1}(y_{sg}),h^{-1}(\Psi_s(y_g)))\in A$ for all $g\in G$, $s\in S$ and hence, $(h^{-1}(y_{sg}),\Phi_s(h^{-1}(y_g)))\in A$ for all $g\in G$, $s\in S$. Thus, $\lbrace h^{-1}(y_g)\rbrace_{g\in G}$ is $A$-pseudo orbit for $\Phi$. 
\medskip

Since $\Phi$ has shadowing, there exists $y_0\in X$ such that $(h^{-1}(y_g),\Phi_g(y_0))\in D$ for all $g\in G$ which implies $(h(h^{-1}(y_g)),h(\Phi_g(y_0)))\in E$ for all $g\in G$. Thus, $(y_g,\Psi_g(h(y_0)))\in E$ for all $g\in G$. This means that $h(y_0)$, $E$-shadows the $D$-pseudo orbit for $\Psi$. This completes the proof. 
\medskip

(iii) Suppose that $\Phi\in Act(G,X)$ is persistent. Let $h:X\rightarrow Y$ be a uniform equivalence satisfying $h\circ\Phi_g=\Psi_g\circ h$ for all $g\in G$. We want to show that $\Psi$ is persistent. Let $E\subset Y\times Y$ be given entourage and let the entourage $D\subset X\times X$ be given by the uniform continuity of $h$, i.e; $(h(x),h(y))\in E$, whenever $(x,y)\in D$. Let the entourage $A\subset X\times X$ be given for $D$ by the persistence of $\Phi$. Since $h^{-1}:Y\rightarrow X$ is uniform equivalence, there exists another entourage $B\subset Y\times Y$ such that $(h^{-1}(x),h^{-1}(y))\in A$, whenever $(x,y)\in B$.  
\medskip

Let $\Psi'\in Act(G,Y)$ be an action such that $(\Psi'_s(y),\Psi_s(y))\in B$ for all $y\in Y$, $s\in S$. For each $s\in S$, let $\Phi'_s=h^{-1}\circ\Psi'_s\circ h$ and consider that $\Phi$ is the action generated by $\lbrace\Phi'_s\mid s\in S\rbrace$. Since $(\Psi'_s(y),\Psi_s(y))\in B$ for all $y\in Y$, $s\in S$ we have $(h^{-1}(\Psi'_s(y)),h^{-1}(\Psi_s(y)))\in A$ for all $y\in Y$, $s\in S$ which implies $(\Phi'_s(h^{-1}(y)),\Phi_s(h^{-1}(y)))\in A$ for all $y\in Y$, $s\in S$. If we put $h(x)=y$, then $(\Phi'_s(x),\Phi_s(x))\in A$ for all $x\in X$, $s\in S$. So by persistence of $\Phi$, for each $x\in X$ there exists $x'\in X$ such that $(\Phi'_g(x),\Phi_g(x'))\in D$ for all $g\in G$. Thus, $(h(\Phi'_g(x)),h(\Phi_g(x')))\in E$ for all $g\in G$ which implies $(\Psi'_g(h(x)),\Psi_g(h(x')))\in E$ for all $g\in G$. This shows that $\Psi$ is persistent.
\medskip

Given a Borel measure $\mu$ on $X$ and a uniform equivalence $h:X\rightarrow Y$, we denote by $h^*(\mu)$ the pullback of $\mu$ defined by $h^*(\mu)(A)=\mu(h^{-1}(A))$ for all Borel measurable set $A\subset Y$. 
\medskip

\textbf{Proposition 4.2} Let $\mu$ be a Borel measure on $X$. Let $\Phi\in Act(G,X)$ and $\Psi\in Act(G,Y)$ be uniformly conjugate. Then,
\begin{enumerate}
\item[(i)] If $\Phi$ has $\mu$-shadowing, then $\Psi$ has $h^*(\mu)$-shadowing. 
\item[(ii)] If $\Phi$ is $\mu$-expansive, then $\Psi$ is $h^*(\mu)$-expansive.
\item[(iii)] If $\Phi$ is $\mu$-persistent, then $\Psi$ is $h^*(\mu)$-persistent. 
\item[(iv)] If $\Phi$ is $\mu$-topologically stable, then $\Psi$ is $h^*(\mu)$-topologically stable.     
\end{enumerate}
\medskip

\textbf{Proof.} We prove only (ii) and (iv). The proofs of (i) and (iii) are left for the reader as easy exercises.   
\medskip

(ii) Suppose that $\Phi$ is $\mu$-expansive. Let $h:X\rightarrow Y$ be a uniform equivalence satisfying $h\circ\Phi_g=\Psi_g\circ h$ for all $g\in G$. We want to show that $\Psi$ is $h^*(\mu)$-expansive. Let $E\subset X\times X$ be a $\mu$-expansive entourage for $\Phi$. Since $h$ is uniform equivalence, there is an entourage $D\subset Y\times Y$ such that $(x,y)\in D$ implies $(h^{-1}(x),h^{-1}(y))\in E$. We claim that $\Gamma_D(y)\subset h(\Gamma_E(h^{-1}(y)))$ for all $y\in Y$. Let $x\in \Gamma_D(y)$, i.e; $(\Psi_g(x),\Psi_g(y))\in D$ for all $g\in G$ which implies $(h^{-1}(\Psi_g(x)),h^{-1}(\Psi_g(y))\in E$ for all $g\in G$. Thus, $(\Phi_g(h^{-1}(x)),\Phi_g(h^{-1}(y)))\in E$ for all $g\in G$ which gives $x\in h(\Gamma_E(h^{-1}(y)))$. So, the claim holds. Therefore, $h^*(\mu)(\Gamma_D(y))\leq \mu(\Gamma_E(h^{-1}(y))$ for all $y\in Y$. So, $h^*(\mu)(\Gamma_D(y))=0$ for all $y\in Y$. Thus, $\Psi$ is $h^*(\mu)$-expansive with $h^*(\mu)$-expansive entourage $D$.  
\medskip

(iv) Suppose that $\Phi$ is $\mu$-topologically stable. Let $h:X\rightarrow Y$ be a uniform equivalence such that $h\circ\Phi_g=\Psi_g\circ h$ for all $g\in G$. We want to show that $\Psi$ is $h^*(\mu)$-topologically stable. Let $E\subset Y\times Y$ be given entourage and let the entourage $D\subset X\times X$ be given by the uniform continuity of $h$, i.e; $(h(x),h(y))\in E$ whenever $(x,y)\in D$. Let the entourage $A\subset X\times X$ be given for $D$ by the $\mu$-topological stability of $\Phi$. Since $h^{-1}:Y\rightarrow X$ is uniform equivalence, there exists another entourage $B\subset Y\times Y$ such that $(h^{-1}(x),h^{-1}(y))\in A$, whenever $(x,y)\in B$.   
\medskip

Let $\Psi'\in Act(G,Y)$ be an action such that $(\Psi'_s(y),\Psi_s(y))\in B$ for all $y\in Y$, $s\in S$. For each $s\in S$, let $\Phi'_s=h^{-1}\circ\Psi'_s\circ h$ and consider that $\Phi$ is the action generated by $\lbrace\Phi'_s\mid s\in S\rbrace$. Since $(\Psi'_s(y),\Psi_s(y))\in B$ for all $y\in Y$, $s\in S$ we have $(h^{-1}(\Psi'_s(y)),h^{-1}(\Psi_s(y)))\in A$ for all $y\in Y$, $s\in S$ which implies $(\Phi'_s(h^{-1}(y)),\Phi_s(h^{-1}(y)))\in A$ for all $y\in Y$, $s\in S$. If we put $h(x)=y$, then $(\Phi'_s(x),\Phi_s(x))\in A$ for all $x\in X$, $s\in S$. So by $\mu$-topological stability of $\Phi$ there exists a compact-valued upper semi-continuous map $H:X\rightarrow \mathcal{P}(X)$ such that $\mu(X\setminus Dom(H))=0$, $\Phi_g\circ H=H\circ\Phi'_g$ for all $g\in G$, $(Id,H)\in D$ and $\mu\circ H=0$.  
\medskip

Put $K=h\circ H\circ h^{-1}$. It is clear that $K$ is upper semi-continuous compact-valued map of $Y$. Since $x\in Dom(K)\Leftrightarrow K(x)\neq\phi\Leftrightarrow h(H(h^{-1}(x)))\neq\phi\Leftrightarrow H(h^{-1}(x))\neq\phi\Leftrightarrow h^{-1}(x)\in Dom(H)\Leftrightarrow x\in h(Dom(H))$. Thus, we have that $Dom(K)=h(Dom(H))$ is measurable. In particular, $K$ has measurable domain. In addition, $h^*(\mu)(Y\setminus Dom(K))=\mu(h^{-1}(Y\setminus Dom(K))=\mu(X\setminus Dom(H))=0$. Observe that $h^*(\mu)(K)(x)=\mu(h^{-1}(K(x)))=\mu(H(h^{-1}(x)))=0$ for all $x\in X$ proving $h^*\mu\circ K=0$. Since $(x,H(x))\in D$ for all $x\in X$, we have $(x,h^{-1}(K(h(x)))\in D$ for all $x\in X$ which implies $(h(x),K(h(x)))\in E$ for all $x\in X$, i.e; $(y,K(y))\in E$ for all $y\in Y$. Finally, $K\circ \Psi'_g=h\circ H\circ h^{-1}\circ\Psi'_g=h\circ H\circ\Phi'_g\circ h^{-1}=h\circ\Phi_g\circ H\circ h^{-1}=\Psi_g\circ h\circ H\circ h^{-1}=\Psi_g\circ K$ for all $g\in G$. This completes the proof.   
\medskip

\textbf{Acknowledgements:} The first author is supported by Department of Science and Technology, Government of India, under INSPIRE Fellowship (Resgistration No-IF150210) program since March 2015.              
       
\end{document}